\documentclass[12pt]{article}
\usepackage{amsmath}\allowdisplaybreaks[4]
\usepackage{amssymb}
\usepackage{amsthm}
\usepackage{amsfonts}
\usepackage{amsmath, dsfont}
\usepackage[margin=1in]{geometry}
\usepackage{graphicx}
\usepackage{graphics}
\usepackage{color}

\newtheorem{theorem}{Theorem}
\newtheorem{corollary}{Corollary}
\begin{document}

\title{Consistent Minimal Displacement of Branching Random Walks}
\author{Ming Fang\thanks{School of Mathematics, University of Minnesota,
206 Church St. SE, Minneapolis, MN 55455, USA.
The work of this author was partially
supported by NSF grant DMS-0804133}
\and
Ofer Zeitouni\thanks{School of Mathematics, University of Minnesota,
206 Church St. SE, Minneapolis, MN 55455, USA and Faculty of Mathematics,
Weizmann
Institute, POB 26, Rehovot 76100, Israel.
The work of this author was partially
supported by NSF grant DMS-0804133}}
\date{December 4, 2009}
\maketitle

\abstract

Let $\mathbb{T}$ denote a rooted $b$-ary tree and let $\{S_v\}_{v\in \mathbb{T}}$ denote a branching random walk indexed by the vertices of the tree, where the increments are i.i.d. and possess a logarithmic moment generating function $\Lambda(\cdot)$. Let $m_n$ denote the minimum of the variables $S_v$ over all vertices at the $n$th generation, denoted by $\mathbb{D}_n$. Under mild conditions, $m_n/n$ converges almost surely to a constant, which for convenience may be taken to be $0$. With $\bar S_v=\max\{S_w:\mbox{\rm $w$ is on the geodesic connecting the root
to $v$}\}$, define $L_n=\min_{v\in \mathbb{D}_n} \bar S_v$. We prove that $L_n/n^{1/3}$ converges almost surely to an explicit constant $l_0$. This answers a question of Hu and Shi.

\section{Introduction}

A branching random walk, as its name suggests, is a process describing a particle performing random walk while branching. In this paper, we consider the 1-dimensional case as follows. At time $0$, there is one particle at location $0$. At time $1$, the particle splits into $b$ particles ($b\in \mathds{Z}_+$ deterministic and $b\geq 2$ to avoid trivial cases), each of which moves independently to a new position according to some distribution function $F(x)$. Then at time 2, each of the $b$ particles splits again into $b$ particles, which again move independently according to the distribution function $F(x)$. The splitting and moving continue at each integer time and are independent of each other. This procedure produces a 1-dimensional branching random walk.

To describe the relation between particles, we associate to each particle a
vertex in a $b$-ary rooted tree $\mathbb{T}=\{V,E\}$ with root $o$,
where each vertex has $b$ children; $V$ is
the set of vertices in $\mathbb{T}$ and $E$ is the set of
edges in $\mathbb{T}$. The root $o$ is associated with the original particle. The $b$ children of a vertex $v\in V$ correspond to the $b$ particles from the splitting of the particle corresponding to $v$. In particular,
the vertices whose distance from $o$ is $n$, denoted by $\mathbb{D}_n$,
correspond to
particles at time $n$. To describe
the displacement between particles, we assign i.i.d. random variables $X_e$
with common
distribution $F(x)$ to each edge $e\in E$. (Throughout, we let $e=uv$ denote the edge $e$ connecting two vertices $u,v\in V$.) For each vertex $v\in V$, we use $|v|$ to denote its distance from $o$ and use $v^k$ to denote the ancestor of $v$ in $\mathbb{D}_k$ for any $0\leq k\leq |v|$. Then the positions of particles at time $n$ can be described by $\{S_v|v\in \mathbb{D}_n\}$, where for $v\in \mathbb{D}_n$, $S_v=\sum_{i=0}^{n-1}X_{v^iv^{i+1}}$.

The limiting behavior of the maximal displacement
$M_n=\max_{v\in \mathbb{D}_n}S_v$ or the minimal
displacement $m_n=\min_{v\in \mathbb{D}_n}S_v$ as $n\to \infty$
has been extensively studied in the literature
(See in particular Bramson \cite{Br1},\cite{Br2},
Addario-Berry and Reed \cite{AR}, and references therein.)
 Throughout this paper, we assume that
\begin{equation}\label{LDA}
Ee^{\lambda X_e}<\infty~~\text{for some}~\lambda<0~\text{and some}~\lambda>0.
\end{equation}
Then the Fenchel-Legendre transform of the log-moment generating
function $\Lambda(\lambda)=\log Ee^{\lambda X_e}$,
\begin{equation}\label{FLT}
\Lambda^*(x)=\sup_{\lambda\in \mathds{R}}(\lambda x-\Lambda(\lambda)),
\end{equation}
is the large deviation rate function
(see \cite[Ch. 1,2]{DZ}) of a
random walk with step distribution $F(x)$.
In addition to \eqref{LDA},
we also assume that, for some $\lambda_-<0$ and
$\lambda_+>0$ in
the interior of $\{\lambda: \Lambda(\lambda)<\infty\}$,
\begin{equation}\label{MaxMinCondition}
 \lambda_{\pm}\Lambda'(\lambda_{\pm})-\Lambda(\lambda_{\pm})=\log b,
\end{equation}
which implies that
$\Lambda^*(\Lambda'(\lambda_{\pm}))=\log b$. These assumptions
imply that
\begin{equation} \label{MaxMin}
  M:=\lim_{n\to \infty}\frac{M_n}{n}=\Lambda'(\lambda_+)~~\text{and}~~
  m:=\lim_{n\to \infty}\frac{m_n}{n}=\Lambda'(\lambda_-)~~~a.s.\:.
\end{equation}
See  \cite{AR}  for
more details on (\ref{MaxMin}).

The \emph{offset} of the branching random walk is defined as the minimal
 deviation of the path up to time $n$ from the line leading to $mn$
(roughly, the minimal position at time n). Explicitly, set
\begin{equation}\label{FD}
L_n=\min_{v\in \mathbb{D}_n}\max_{k=0}^{n}(S_{v^k}-mk).
\end{equation}
Without loss of generality, subtracting the deterministic constant
$\Lambda'(\lambda_-)$ from
 each increment $\{X_e\}$, we can and will assume that
\begin{equation}\label{ZeroMin}
m=\Lambda'(\lambda_-)=0.
\end{equation}
Under this assumption, (\ref{MaxMinCondition}) and (\ref{FD}) simplify to
$$-\Lambda(\lambda_-)=\log b,~\eqno{(\ref{MaxMinCondition}')}$$
$$L_n=\min_{v\in \mathbb{D}_n}\max_{k=0}^{n}S_{v^k}.~\eqno{(\ref{FD}')}$$
In the process of studying random walks in random environments on trees,
Hu and Shi \cite{HS1} (2007) discovered that the offset has order $n^{1/3}$ in the following sense: there exist constants $c_1,c_2>0$ such that
\begin{equation}\label{HS}
c_1\leq \liminf_{n\to \infty}\frac{L_n}{n^{1/3}}\leq \limsup_{n\to \infty}\frac{L_n}{n^{1/3}}\leq c_2.
\end{equation}
They raised and advertised the question as to whether the limit of
${L_n}/{n^{1/3}}$ exists. In this note, we answer this affirmatively
and prove the following.
\begin{theorem}
\label{theo-main}
Under assumption (\ref{LDA}) and (\ref{MaxMinCondition}) and with
$l_0=\sqrt[3]{\frac{3\pi^2\sigma_Q^2}{-2\lambda_-}}$, it holds that
\begin{equation}\label{THM}
  \lim_{n\to \infty} \frac{L_n}{n^{1/3}}=l_0~~~a.s.\;.
\end{equation}
\end{theorem}
In the expression for  $l_0$,
$\lambda_-<0$ by the definition (\ref{MaxMinCondition})
 and $\sigma_Q^2$ is a certain variance defined in (\ref{variance}).

The proof of the theorem is divided into two parts - the lower bound
(\ref{lb}) and the upper bound (\ref{ub}). In Section 2, we review a
result from Mogul'skii \cite{Mo}, which will be the key estimate in
our proof. In Section 3, we apply a first moment argument (with a twist)
in order
to
study
the minimal positions for intermediate levels with the
restriction that the walks do not exceed $ln^{1/3}$ for some $l>0$ at
all time. This yields the lower bound for $L_n$.
In section 4, we apply a second moment argument to lower bound
 $P(L_n\leq l n^{1/3})$ for certain values of $l$.
Compared with standard applications of the second moment method in
related problems, the analysis here requires the control of second order
terms in the large deviation estimates. Truncation of the tree is then used
to get independence and complete the proof of the upper bound.

\section*{Acknowledgement}
While this work was being
completed, we learnt
that
as part of their study of RWRE on trees,
G. Faraud, Y. Hu  and Z. Shi had independently obtained
Theorem \ref{theo-main},
using a related but slightly different method \cite{FHZ}. In particular,
their work handles also the case of Galton--Watson trees.
We thank Y. Hu for discussing this problem with one of us (O.Z.) and for
providing us with the reference \cite{Mo}, which allowed
us to skip tedious details in our original proof.

\section{An Auxiliary Estimate: the absorption problem for random walk}

We derive in this section some estimates for random walk with i.i.d.
increments
$\{X_i\}_{i\geq 1}$
distributed according to a law $P$ with $P((-\infty,x])=F(x)$ satisfying
\eqref{LDA}, \eqref{MaxMinCondition} and \eqref{ZeroMin}.
Define
$$S_n(t)=\frac{X_0+X_1+\dots+X_k}{n^{1/3}}~~~\\
\text{for}~\frac{k}{n}\leq t<\frac{k+1}{n},~ k=0,1,\dots,n-1,$$
where $X_0=0$. Note that due to (\ref{ZeroMin}), $EX_i>0$.
 Introduce the auxiliary law
\begin{equation}\label{CoM}
  \frac{d Q}{d P}=e^{\lambda_-X_1-\Lambda(\lambda_-)}.
\end{equation}
Under $Q$, $E_QX_1=0$. The variance of $X_1$ under $Q$ is denoted by
\begin{equation}\label{variance}
   \sigma_Q^2=E_QX_1^2.
\end{equation}
In the following estimates, $f_1(t)$ and $f_2(t)$, which may take the value $\pm\infty$, are right-continuous and piecewise constant functions on $[0,1]$. $G=\cup_{0\leq t\leq 1}\{(f_1(t),f_2(t))\times t\}$ is a region bounded by $f_1(t)$ and $f_2(t)$. Assume also that $G$ contains the graph of a continuous function.
\begin{theorem}(Mogul'skii \cite[Theorem 3]{Mo})
\label{theo-mog}
Under the above assumptions,
\begin{equation}\label{Mogulskii}
Q(S_n(t)\in G, t\in[0,1])=
e^{-\frac{\pi^2\sigma_Q^2}{2}H_2(G)n^{1/3}+o(n^{1/3})},
\end{equation}
where
\begin{equation}\label{H2G}
  H_2(G)=\int_0^1\frac{1}{(f_1(t)-f_2(t))^2}dt.
\end{equation}
\end{theorem}
In the following, we will need to control the dependence of the estimate
\eqref{Mogulskii} on the starting point.
\begin{corollary}
With notation and assumptions as in Theorem
\ref{theo-mog},
for any $\epsilon>0$, there is a $\delta>0$ such that,
for any interval $I\subset (f_1(0),f_2(0))$ with length $|I|\leq \delta$,
we have
\begin{equation}\label{supbound}
  \sup_{x\in I}Q(x+S_n(\cdot)\in G)\leq
e^{-(\frac{\pi^2\sigma_Q^2}{2}H_2(G)-\epsilon)n^{1/3}+o(n^{1/3})}
.
\end{equation}
\end{corollary}
\noindent
{\bf Proof}  Let $I=(a,b)$ and
$G_x:=\cup_{0\leq t\leq 1}\{(f_1(t)-x,f_2(t)-x)\times t\}$ be the
 shift of $G$ by $x$. Set $G'=G_a\cup G_b$. We have
$$\sup_{x\in I}Q(x+S_n(\cdot)\in G)=\sup_{x\in I} Q(S_n(\cdot)\in G_x)
\leq
Q(S_n(\cdot)\in G')=e^{-\frac{\pi^2\sigma_Q^2}{2}H_2(G')n^{1/3}+o(n^{1/3})}
\,.$$
\\
Since $H_2(G')=\int_0^1\frac{1}{(f_2(t)-f_1(t)+(b-a))^2}dt\uparrow H_2(G)$ as $|I|=(b-a)\to 0$ uniformly in the position of $I$, the lemma is proved.
\qed

\section{Lower Bound}

Consider the branching random walk up to level $n$. In this and the next
section, we estimate the number of particles that stay
constantly below $ln^{1/3}$, i.e.,
\begin{equation}\label{N1}
N^l_n=\sum_{v\in \mathbb{D}_n}
1_{\{S_{v^k}\leq ln^{1/3}~\text{for}~k=0,1,\dots,n\}}.
\end{equation}
In order to get a lower bound on the offset, we apply a first moment
method with a small twist: while it is natural to just calculate the
first moment of $N_n^l$, such a computation ignores the constraint on
 the number of particles at level $k$ imposed by the tree structure.
In particular, $EN^l_n$ for branching random walks is the same as the one
 for $b^n$ independent random walks. An easy first and second moment
 argument shows that the limit in (\ref{THM}) is 0 for $b^n$ independent
 random walks, and thus no useful upper bound can be derived in this way.

To address this issue,
we use a more delicate first moment argument. Namely, we look at the
 vertices not only at level $n$ but also at some intermediate levels.
Divide the interval $[0,n]$ into $1/\epsilon$ equidistant levels,
with $1/\epsilon$ an integer. Define recursively, for any $\delta>0$,
\begin{equation}\label{s w}
 \left\{\begin{array}{l}
          s_0=0,w_0=l+\delta;\\
          s_k=s_{k-1}-\frac{\pi^2\sigma_Q^2}{2\lambda_-w_{k-1}^2}\epsilon,~w_k=l+\delta-s_k~\text{for}
          ~k=1,\dots,\frac{1}{\epsilon}.
        \end{array}
 \right.
\end{equation}
For particles staying
below $ln^{1/3}$,
$s_k$ will be interpreted
as
values such that the walks between times
$k\epsilon n$ and $(k+1)\epsilon n$ never go below $(s_k-\delta) n^{1/3}$,
and $w_kn^{1/3}$ will correspond to
the width of the
window $W_k=((s_k-\delta) n^{1/3},ln^{1/3})$ that
we allow between level $k\epsilon n$ and $(k+1)\epsilon n$,
when considering those
walks that do not go below $(s_k-\delta) n^{1/3}$ or go above $ln^{1/3}$.

Before calculating the first moment,
consider the recursion (\ref{s w}) for $s_k$. Rewrite it as
\begin{equation}\label{recursion s}
 s_k=s_{k-1}-\frac{\pi^2\sigma_Q^2}{2\lambda_-(l+\delta -s_{k-1})^2}\epsilon.
\end{equation}
This is an Euler's approximation sequence for the solution of the
following differential equation
\begin{equation}\label{DE}
s'(t)=-\frac{\pi^2\sigma_Q^2}{2\lambda_-(\alpha -s(t))^2},~s(0)=0,
\end{equation}
where $\alpha=l+\delta$. The above initial value problem has the
 solution $s^{\alpha}(t)=\alpha+\sqrt[3]{-\frac{3\pi^2\sigma_Q^2}{2\lambda_-}t-\alpha^3}$. Here we find
\begin{equation}\label{l0}
  l_0=\sqrt[3]{\frac{3\pi^2\sigma_Q^2}{-2\lambda_-}}
\end{equation}
such that $s^{l_0}(1)=l_0$. For any $\l_1<l_0$, we can choose $\delta>0$ and $l_1+\delta<l_0$. In this case, $s^{l_1+\delta}(1)>l_1+\delta>l_1$. If we choose such $l_1$ and $\delta$ in (\ref{s w}), it is easy to check that the sequence $\{s_k\}_{k=0}^{\frac{1}{\epsilon}}$ will be greater than $l_1$ somewhere in the sequence. Define
\begin{equation}\label{KDef}
  K=\min\{k:s_k\geq l_1\}.
\end{equation}
For fixed $\gamma>0$ small enough, we can choose $\epsilon$ small such that
\begin{equation}\label{gamma}
  K\epsilon<1-\gamma.
\end{equation}
For $k<K-1$,
let $Z_k$ denote the number of vertices $v$
between level
$k\epsilon n$ and $(k+1)\epsilon n$ with
$S_v<(s_k-\delta)n^{1/3}$.
Denote by $Z_{K-1}$ the number of vertices $w$
between level $(K-1)\epsilon n$ and $n$ with
$S_w<(s_{K-1}-\delta)n^{1/3}$. Denote by $Z$ the number
vertices $v\in \mathbb{D}_n$
whose associated
walks stay in $W_k$ between level $k\epsilon n$ and $(k+1)\epsilon n$
for $k<K$ and then stay in $W_{K-1}$ up to level $n$. Explicitly,

\begin{eqnarray}
      &&Z_0=\sum_{i=1}^{\lfloor
\epsilon n\rfloor}\sum_{v\in \mathbb{D}_i}1_{\{S_v<-\delta n^{1/3}\}},\\
      &&Z_k=\sum_{i=\lfloor
k\epsilon n\rfloor
+1}^{\lfloor
(k+1)\epsilon n\rfloor}
\sum_{v\in \mathbb{D}_i}1_{\{S_v<(s_k-\delta)n^{1/3},~S_{v^d}\in W_{j}~\text{for}~j\epsilon n \leq d\leq (j+1)\epsilon n ~ \text{and}~ j<k\}},~~~0<k<K-1,\\
      &&Z_{K-1}=\sum_{i=\lfloor K\epsilon n\rfloor
+1}^{n}\sum_{v\in \mathbb{D}_i}1_{\{S_v<(s_{K-1}-\delta)n^{1/3},~S_{v^d}\in W_{j}~\text{for}~j\epsilon n \leq d\leq (j+1)\epsilon n ~ \text{and}~ j<K-1\}},\\
      &&Z=\sum_{v\in\mathbb{D}_n}1_{\{S_{v^d}\in W_{j}~\text{for}~j\epsilon n \leq d\leq (j+1)\epsilon n ~ \text{and}~ j<K,~S_{v^d}\in W_{K-1}~\text{for}~K\epsilon n \leq d\leq n \}}.
\end{eqnarray}
\\
\begin{figure}
   \centering
   \includegraphics[width=90mm]{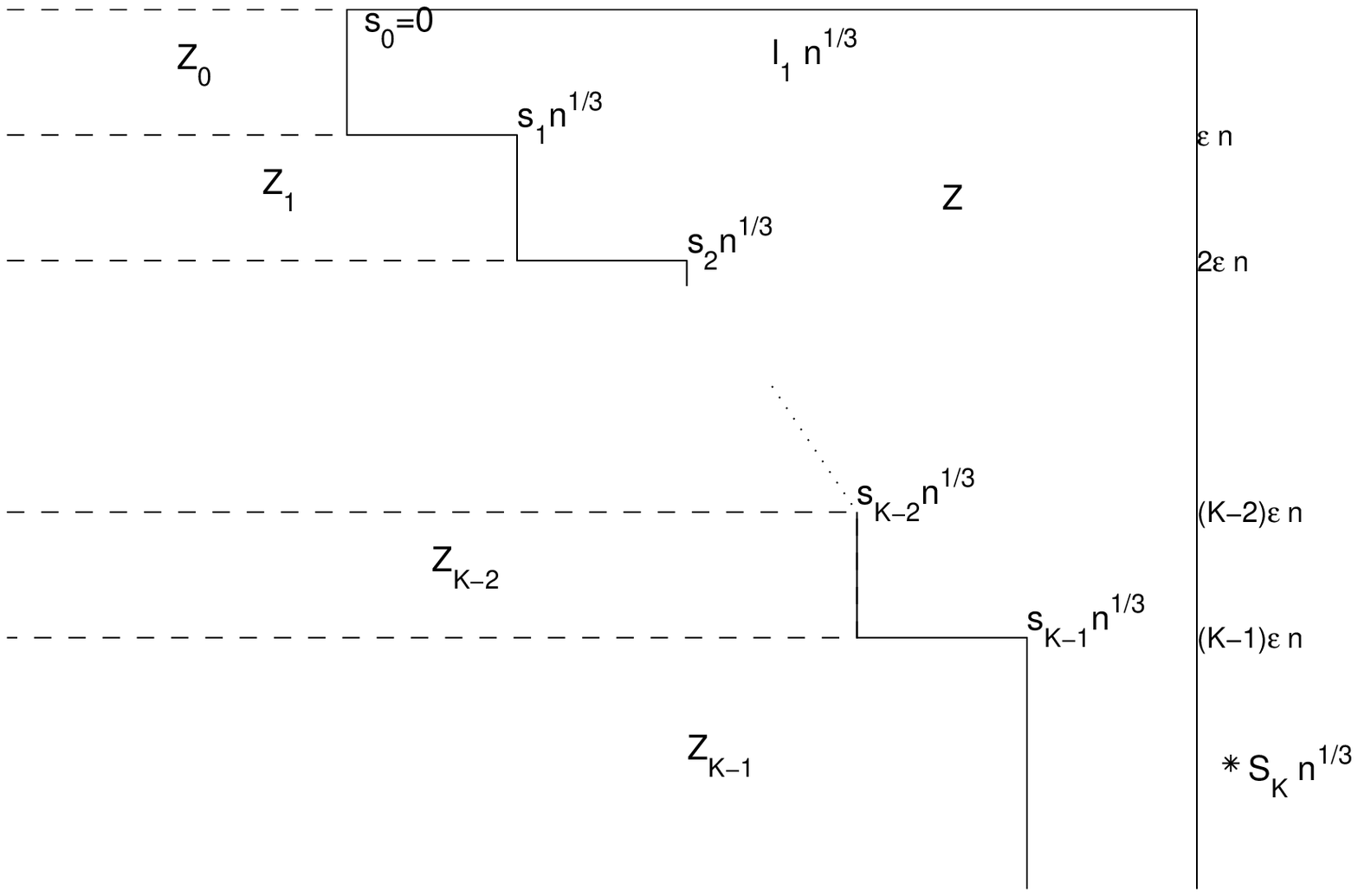}
   \caption{The relation between $Z_k$'s and $s_k$'s.}
\end{figure}
\\
Observe that $N^{l_1}_n\leq \sum_{k=0}^{K-1}Z_k+Z$.
Using Theorem \ref{theo-mog},
we provide upper bounds for the first moment of the
$Z_k$s and $Z$. Starting with $Z_0$, we have
\begin{eqnarray*}
   EZ_0&=&\sum_{i=1}^{\lfloor
\epsilon n\rfloor} b^i E1_{\{S_i<-\delta n^{1/3}\}}= \sum_{i=1}^{\lfloor
\epsilon n\rfloor} b^i E_Q e^{-\lambda_-S_i+i\Lambda(\lambda_-)}1_{\{S_i<-\delta n^{1/3}\}}\\
   &\leq&\sum_{i=1}^{\lfloor
\epsilon n\rfloor} e^{\lambda_-\delta n^{1/3}}E_Q
1_{\{S_i<-\delta n^{1/3}\}}\leq
\sum_{i=1}^{\lfloor\epsilon n\rfloor} e^{\lambda_-\delta n^{1/3}}
   \leq  e^{\lambda_-\delta n^{1/3}+o(n^{1/3})},
\end{eqnarray*}
where we used the change of measure (\ref{CoM}) in the second equality, and ($3'$) and the fact that $\lambda_-<0$ in the first inequality. For $0<k<K-1$, using again the change of measure (\ref{CoM}), we get
\begin{eqnarray*}
 E Z_k &=& \sum_{i=\lfloor
k\epsilon n\rfloor
+1}^{\lfloor
(k+1)\epsilon n\rfloor
} b^{i} E 1_{\{S_i<(s_k-\delta)n^{1/3},~S_d\in W_{j}~\text{for}~j\epsilon n \leq d\leq (j+1)\epsilon n ~ \text{and}~ j<k\}}\\
 &=& \sum_{i=\lfloor
k\epsilon n\rfloor
+1}^{\lfloor
(k+1)\epsilon n\rfloor
} E_Q e^{-\lambda_-S_i}1_{\{S_i<(s_k-\delta)n^{1/3},~S_d\in W_{j}~\text{for}~j\epsilon n \leq d\leq (j+1)\epsilon n ~ \text{and}~ j<k\}}\\
 &\leq & e^{-\lambda_-(s_k-\delta)n^{1/3}}\sum_{i=\lfloor
k\epsilon n\rfloor
+1}^{\lfloor
(k+1)\epsilon n\rfloor
}E_Q 1_{\{S_i<(s_k-\delta)n^{1/3},~S_d\in W_{j}~\text{for}~j\epsilon n \leq d\leq (j+1)\epsilon n ~ \text{and}~ j<k\}}.
\end{eqnarray*}
Therefore,
\begin{eqnarray*}
 EZ_k&\leq& e^{-\lambda_-(s_k-\delta)n^{1/3}}\sum_{i=\lfloor
k\epsilon n\rfloor+1}^{\lfloor(k+1)\epsilon ni\rfloor
} Q(S_d\in W_{j}~\text{for}~j\epsilon n \leq d\leq (j+1)\epsilon n ~ \text{and}~ j<k)\\
 &=& e^{-\lambda_-(s_k-\delta)n^{1/3}}\sum_{i=\lfloor
k\epsilon n\rfloor
+1}^{\lfloor
(k+1)\epsilon n\rfloor} e^{-\sum_{j=0}^{k-1}\frac{\pi^2\sigma_Q^2}{2w_j^2}\epsilon n^{1/3}+o(n^{1/3})}\\
 &\leq & e^{-\lambda_-(s_k-\delta)n^{1/3}-\sum_{j=0}^{k-1}\frac{\pi^2\sigma_Q^2}{2w_j^2}\epsilon n^{1/3}+o(n^{1/3})}
 =e^{\lambda_-\delta n^{1/3}+o(n^{1/3})},
\end{eqnarray*}
where (\ref{Mogulskii}) with the choice of $G=\{\cup_{j=0}^{k-1}W_j/n^{1/3}\times [j\epsilon ,(j+1)\epsilon )\}\cup \{(-\infty,\infty)\times[k\epsilon,1]\}$ is applied in the first equality, and (\ref{s w}) in the second. The calculation of $EZ_{K-1}$ is almost the same as $EZ_k$ except that we replace the summation limits above by $(K-1)\epsilon n+1$ and $n$ and that we replace the $k$ in the summand by $K-1$. Thus, we get the same upper bound for $EZ_{K-1}$,
$$EZ_{K-1}\leq e^{\lambda_-\delta n^{1/3}+o(n^{1/3})}.$$
We estimate $EZ$ similarly as follows. First, use the change
of measure
(\ref{CoM}) to get
%
%
%
\begin{eqnarray*}
   EZ & = & b^nE1_{\{S_d\in W_{j}~\text{for}~j\epsilon n \leq d\leq (j+1)\epsilon n ~ \text{and}~ j<K,~S_d\in W_{K-1},~\text{for}~K\epsilon n \leq d\leq n \}}\\
   &=& E_Q e^{-\lambda_-S_n} 1_{\{S_d\in W_{j}~\text{for}~j\epsilon n
\leq d\leq (j+1)\epsilon n ~ \text{and}~ j<K,~S_d\in W_{K-1},~\text{for}~K
\epsilon n \leq d\leq n \}}\\
   &\leq & e^{-\lambda_-l_1n^{1/3}}E_Q 1_{\{S_d\in W_{j}~\text{for}~j
\epsilon n \leq d\leq (j+1)\epsilon n ~ \text{and}~ j<K,~S_d\in
W_{K-1},~\text{for}~K\epsilon n \leq d\leq n \}}.
\end{eqnarray*}
Then, applying
(\ref{Mogulskii}) with $G=\{\cup_{j=0}^{K-1}W_j/n^{1/3}\times [j\epsilon ,(j+1)\epsilon )\}\cup \{W_{K-1}/n^{1/3}\times[K\epsilon,1]\}$ in the first
equality, we get
\begin{eqnarray*}
   EZ & \leq &
    e^{-\lambda_-l_1n^{1/3}}E_Q 1_{\{S_d\in W_{j}~\text{for}~j\epsilon n \leq d\leq (j+1)\epsilon n ~ \text{and}~ j<K,~S_d\in W_{K-1},~\text{for}~K\epsilon n \leq d\leq n \}}\\
   & = & e^{-\lambda_-l_1n^{1/3}-\sum_{i=0}^{K-1}\frac{\pi^2\sigma_Q^2}{2w_i^2}\epsilon n^{1/3}-\frac{\pi^2\sigma_Q^2}{2w_{K-1}^2}(1-K\epsilon)n^{1/3}+o(n^{1/3})}\\
   &\leq & e^{-\gamma\frac{\pi^2\sigma_Q^2}{2l_1^2}n^{1/3}+o(n^{1/3})},
\end{eqnarray*}
where the last inequality is obtained by noting that
$l_1\leq S_K=-\sum_{i=0}^{K-1}\frac{\pi^2\sigma_Q^2}{2\lambda_-w_i^2}\epsilon$
 by (\ref{KDef}) and (\ref{s w}), and then recalling
 (\ref{gamma}) and $w_{K-1}<l_1$.

In conclusion,
we proved that $E(\sum_{k=0}^{K-1}Z_k+Z)\leq e^{-c_3n^{1/3}+o(n^{1/3})}$
for some $0<c_3<\min\{-\lambda_-\delta,\gamma\frac{\pi^2\sigma_Q^2}{2l_1^2}\}$. Since $\sum_{k=0}^{K-1}Z_k+Z$ is an integer valued random variable, we have
$$P(\sum_{k=0}^{K-1}Z_k+Z>0)=P(\sum_{k=0}^{K-1}Z_k+Z\geq 1)\leq E(\sum_{k=0}^{K-1}Z_k+Z)\leq e^{-c_3n^{1/3}+o(n^{1/3})}.$$
By the Borel-Cantelli lemma, we have $\sum_{k=0}^{K-1}Z_k+Z=0$ a.s. for all large $n$. So is $N_n^{l_1}=0$, which means that $L_n>l_1n^{1/3}$ a.s. for all large $n$. Since $l_1<l_0$ is arbitrary, we conclude that
\begin{equation}\label{lb}
  \liminf_{n\to \infty}\frac{L_n}{n^{1/3}}\geq l_0~~~~a.s..
\end{equation}
This completes the proof of the lower bound in Theorem \ref{theo-main}.

\section{Upper Bound}

\subsection{A Second Moment Method Estimate}

In this section, we consider any fixed $l_2>l_0$. A second moment argument
will provide
a lower bound for the probability that we can find at least
one walk which stays in the
interval $W_k$ between level $k\epsilon n$
and $(k+1)\epsilon n$ for all $k$. A truncation (of the tree) argument
will complete the proof of the upper bound.

As a first step, consider the sequence $\{s_k\}$ in (\ref{s w}) with $l_2>l_0$. Then for any $\delta>0$, it is easy to see that $s^{l_2+\delta}(t)$ is increasing and convex for $0\leq t\leq 1$. Thus in Euler's approximation,
\begin{equation}\label{sleql2}
  s_{\frac{1}{\epsilon}}<s^{l_2+\delta}(1)<s^{l_2}(1)<l_2.
\end{equation}
It follows from \eqref{s w} that
\begin{equation}\label{wgeqdelta}
  w_k\geq \delta \;\;\:\text{for all}\;\:0\leq k\leq \frac{1}{\epsilon}-1 .
\end{equation}
Define $\tilde{N}_n^{l_2}$ as follows.
$$\tilde{N}_n^{l_2}=\sum_{v\in \mathbb{D}_n}1_{\{S_{v^j}\in W_k,~\text{for}~k\epsilon n\leq j\leq (k+1)\epsilon n,~ k=0,\dots, \frac{1}{\epsilon}-1\}}. $$
\\
We will apply second moment method to $\tilde{N}_n^{l_2}$. $E\tilde{N}_n^{l_2}$ is calculated the same way as $EZ$ in the previous section. But this time we consider $G=\{\cup_{j=0}^{\frac{1}{\epsilon}-1}W_j/n^{1/3}\times [j\epsilon ,(j+1)\epsilon )\}\cup \{(l_2-\Delta l_2,l_2)\times\{1\}\}$ in (\ref{Mogulskii}) with $\Delta l_2\to 0$, so

\begin{eqnarray}\label{FM}
     E\tilde{N}_n^{l_2} &=& b^n E 1_{\{S_j\in W_k,~\text{for}~k\epsilon n\leq j\leq (k+1)\epsilon n,~ k=0,\dots, \frac{1}{\epsilon}-1\}} \nonumber\\
   &=& E_Qe^{-\lambda_-S_n} 1_{\{S_j\in W_k,~\text{for}~k\epsilon n\leq j\leq (k+1)\epsilon n,~ k=0,\dots, \frac{1}{\epsilon}-1\}}\nonumber\\
   &=& e^{(-\lambda_-l_2
   -\sum_{k=0}^{\frac{1}{\epsilon}-1}\frac{\pi^2\sigma_Q^2}{2w_k^2}\epsilon)n^{1/3}+o(n^{1/3})}.
\end{eqnarray}
\\
From (\ref{sleql2}) and the definition (\ref{s w}) of $s_k$, $-\lambda_-l_2-\sum_{k=0}^{\frac{1}{\epsilon}-1}\frac{\pi^2\sigma_Q^2}{2w_k^2}\epsilon>0$ and thus $E\tilde{N}_n^{l_2}\to \infty$. Therefore, we will be ready to apply the second moment method after the following calculations.

\begin{eqnarray}\label{SM0}
   E (\tilde{N}_n^{l_2})^2&=& E \sum_{u,v\in \mathbb{D}_n} 1_{\{S_{u^j},S_{v^j}\in W_k,~\text{for}~k\epsilon n\leq j\leq (k+1)\epsilon n,~ k=0,\dots, \frac{1}{\epsilon}-1\}}\nonumber\\
   &=& \sum_{h=0}^{n-1}E\sum_{\substack{u,v\in \mathbb{D}_n\\ u\wedge v\in \mathbb{D}_h}} 1_{\{S_{u^j},S_{v^j}\in W_k,~\text{for}~k\epsilon n\leq j\leq (k+1)\epsilon n,~ k=0,\dots, \frac{1}{\epsilon}-1\}}+E\tilde{N}_n^{l_2}.
\end{eqnarray}
\\
In the last expression above, $u\wedge v$ is the largest
common ancestor of $u$ and $v$. Write $h=q\epsilon n+r$ for $0\leq q\leq \frac{1}{\epsilon}-1$ and $0\leq r<\epsilon n$.
There are
$b^{2n-h-1}(b-1)$ indices in the second sum in the right side
of (\ref{SM0}). We estimate the probability for one such pair to stay
in $W_k$'s. In order to simplify the notation,
define
$$p_1(0,h,x)=P(S_h\in dx,~S_j\in W_k,~\text{for}~k\epsilon n\leq j\leq (k+1)\epsilon n\wedge h,~ k=0,\dots,q),$$
$$p_2(h,x,n,y)=P(S_n\in dy,~S_j\in W_k,~\text{for}~h\vee k\epsilon n\leq j\leq (k+1)\epsilon n,~ k=q,\dots,n|S_h=x).$$
\\
Similarly, define $q_1(0,h,x)$ and $q_2(h,x,n,y)$ to be the probability of the same events under $Q$. Then we have

\begin{eqnarray}\label{SM}
   E(\tilde{N}_n^{l_2})^2 & = & E\tilde{N}_n^{l_2}+\sum_{h=0}^{n-1}b^{2n-h-1}(b-1)\int_{W_q}
   (\int_{W_n}p_2(h,x,n,y)dy)^2p_1(0,h,x)dx\nonumber\\
   &=& E\tilde{N}_n^{l_2} +\sum_{h=0}^{n-1} b^{2n-h-1}(b-1)\int_{W_q}(\int_{W_n}
   e^{-\lambda_-(y-x)+(n-h)\Lambda(\lambda_-)}q_2(h,x,n,y)dy)^2\nonumber\\
   &&\;\;\;\;\;\;\;\;\;\;\;\;\;\;\;\;\;\;\;\;\;\;\;\;\;\;\;\;\;\;\;\;\;\;\;\;\;\;
   \;\;\;\;\;\;\;\;\;\cdot   e^{-\lambda_-x+h\Lambda(\lambda_-)}q_1(0,h,x)dx\nonumber\\
   &\leq & E\tilde{N}_n^{l_2}+\sum_{h=0}^{n-1}\frac{b-1}{b}
   e^{(-2\lambda_-l_2+\lambda_-(s_q-\delta))n^{1/3}}\nonumber\\
   &&\:\:\:\:\:\:\:\:\:\:\:\:\:\:\:\:\:\:\:\:\:\:\:\:\:\:\:\:\:\:\:\:\:\:\:\:
   \cdot\int_{W_q}(\int_{W_n}q_2(h,x,n,y)dy)^2q_1(0,h,x)dx.
\end{eqnarray}
\\
We now provide
an upper bound for the integral term in the right side of (\ref{SM}). We have
\begin{eqnarray}
  & & \int_{W_q}(\int_{W_n}q_2(h,x,n,y)dy)^2q_1(0,h,x)dx\nonumber\\
  &\leq & (\sup_{x\in W_q}\int_{W_n}q_2(h,x,n,y)dy)^2\int_{W_q}q_1(0,h,x)dx\nonumber\\
  &\leq & (\sup_{x\in W_q}\int_{W_n}\int_{W_{q+1}}q_2(h,x,(q+1)\epsilon n,z)q_2((q+1)\epsilon n,z,n,y)dzdy)^2\int_{W_q}q_1(0,q\epsilon n,x)dx\nonumber\\
  &\leq & (\sup_{x\in W_q}\int_{W_{q+1}}q_2(h,x,(q+1)\epsilon n,z)q_2((q+1)\epsilon n,z,n,W_n)dz)^2e^{-\sum_{k=0}^{q-1}\frac{\pi^2\sigma_Q^2}{2w_k^2}n^{1/3}+o(n^{1/3})}\nonumber\\
  &\leq & (\sup_{x\in W_q}\sum_{i}\int_{I_i}q_2(h,x,(q+1)\epsilon n,z)q_2((q+1)\epsilon n,z,n,W_n)dz)^2e^{-\sum_{k=0}^{q-1}\frac{\pi^2\sigma_Q^2}{2w_k^2}n^{1/3}+o(n^{1/3})}\nonumber\\
  &\leq & (\sum_{i}\sup_{z\in I_i}q_2((q+1)\epsilon n,z,n,W_n))^2e^{-\sum_{k=0}^{q-1}\frac{\pi^2\sigma_Q^2}{2w_k^2}\epsilon n^{1/3}+o(n^{1/3})}.
\end{eqnarray}
In the above, $\cup_{i} I_i=W_{q+1}$. Due to (\ref{supbound}),
for any small $\epsilon_1>0$, we can choose a
finite number of $I_i$s and $|I_i|\leq \delta_1n^{1/3}$ such that
for each $i$,
$$\sup_{z\in I_i}(q_2((q+1)\epsilon n,z,n,W_n)\leq e^{-(\sum_{k=q+1}^{\frac{1}{\epsilon}-1}\frac{\pi^2\sigma_Q^2}{2w_k^2}\epsilon-\epsilon_1)
n^{1/3}+o(n^{1/3})}.$$
After splitting $\sum_{h=0}^{n-1}$ to $\sum_{q=0}^{\frac{1}{\epsilon}-1}\sum_{r=0}^{\epsilon n-1}$ in (\ref{SM}), we obtain the upper bound for $E(\tilde{N}_n^{l_2})^2$ as follows,
\begin{eqnarray}\label{SMUB}
  E(\tilde{N}_n^{l_2})^2 & \leq & E\tilde{N}_n^{l_2}+\sum_{q=0}^{1/\epsilon-1}
e^{(-2\lambda_-l_2+\lambda_-(s_q-\delta))n^{1/3}
-\sum_{k=0}^{q-1}\frac{\pi^2\sigma_Q^2}{2w_k^2}\epsilon n^{1/3}
-2(\sum_{k=q+1}^{\frac{1}{\epsilon}-1}\frac{\pi^2\sigma_Q^2}{2w_k^2}\epsilon-\epsilon_1)n^{1/3}+o(n^{1/3})}\nonumber\\
&\leq & \sum_{q=0}^{1/\epsilon-1} e^{(-2\lambda_-l_2+\lambda_-(s_q-\delta))n^{1/3}
-\sum_{k=0}^{q-1}\frac{\pi^2\sigma_Q^2}{2w_k^2}\epsilon n^{1/3}
-2(\sum_{k=q+1}^{\frac{1}{\epsilon}-1}\frac{\pi^2\sigma_Q^2}{2w_k^2}\epsilon-\epsilon_1)n^{1/3}+o(n^{1/3})}.
\end{eqnarray}
With the bounds for $E\tilde{N}_n^{l_2}$ (\ref{FM}) and $E(\tilde{N}_n^{l_2})^2$ (\ref{SMUB}), we have
\begin{eqnarray}\label{SMM}
P(\tilde{N}_n^{l_2}>0)
&\geq & \frac{(E\tilde{N}_n^{l_2})^2}{E(\tilde{N_n^{l_2}})^2}\geq \frac{1}{\sum_{q=0}^{1/\epsilon-1}
e^{(\lambda_-(s_q-\delta)+\sum_{k=0}^{q-1}\frac{\pi^2\sigma_Q^2}{2w_k^2}\epsilon
+2\frac{\pi^2\sigma_Q^2}{2w_q^2}\epsilon+2\epsilon_1)n^{1/3}+o(n^{1/3})}}\nonumber\\
&=& \frac{1}{\sum_{q=0}^{1/\epsilon -1}e^{(-\lambda_-\delta +\frac{\pi^2\sigma_Q^2}{w_q^2}\epsilon+2\epsilon_1) n^{1/3}+o(n^{1/3})}} \geq e^{(\lambda_-\delta-\frac{\pi^2\sigma_Q^2}{\delta^2}\epsilon-2\epsilon_1)n^{1/3}+o(n^{1/3})}\nonumber\\
& = & e^{-\epsilon_2n^{1/3}+o(n^{1/3})},
\end{eqnarray}
where $\epsilon_2:=-\lambda_-\delta+\frac{\pi^2\sigma_Q^2}{\delta^2}\epsilon+2\epsilon_1$, and we use (\ref{s w}) in the first equality and $w_q\geq \delta$ (see (\ref{wgeqdelta})) in the last inequality. We can make $\epsilon_2$ arbitrarily small by first choosing $\delta$ small then choosing $\epsilon$ and $\epsilon_1$ small. Therefore, we get
\begin{equation}\label{PLB}
P(L_n\leq l_2n^{1/3})\geq P(\tilde{N}_n^{l_2}>0)\geq e^{-\epsilon_2n^{1/3}+o(n^{1/3})}.
\end{equation}
\subsection{A Truncation Argument}

In view of the lower bound (\ref{PLB}), we truncate the tree at level $\lfloor \epsilon_3n^{1/3}\rfloor=\lfloor 2\epsilon_2 n^{1/3}/{\log b}\rfloor$ to get $b^{\lfloor \epsilon_3 n^{1/3}\rfloor}\geq e^{2\epsilon_2 n^{1/3}}/b$ independent branching random walks. We take care of the path before and after level $\lfloor\epsilon_3n^{1/3}\rfloor$ separately.

Define $L_n^v$ similarly as $L_n$ for each branching random walk starting from $v\in \mathbb{D}_{\lfloor\epsilon_3 n^{1/3}\rfloor}$. Then
\begin{eqnarray}
  P(L_n^v> l_2 n^{1/3}~\text{for every}~ v)&=&(1-P(L_n\leq l_2n^{1/3}))^{b^{\lfloor\epsilon_3n^{1/3}\rfloor}}\nonumber\\
  &\leq & (1-e^{-\epsilon_2n^{1/3}+o(n^{1/3})})^{e^{2\epsilon_2n^{1/3}}/b}\leq e^{-e^{\epsilon_2n^{1/3}+o(n^{1/3})}},
\end{eqnarray}
when $n$ is large. By the Borel-Cantelli lemma, the above double exponential guarantees that almost surely for all large n, there exists a $v\in\mathbb{D}_{\lfloor\epsilon_3n^{1/3}\rfloor}$ such that
\begin{equation}\label{LVUB}
  L_n^v\leq l_2n^{1/3}.
\end{equation}
This is an upper bound for the deviation of paths after
level $\lfloor\epsilon_3n^{1/3}\rfloor$. We also need to control the paths before that
level, which is a standard large deviation computation. Indeed, for $q$
integer (later, we take $q=
 \lfloor\epsilon_3n^{1/3}\rfloor$), set
%
$$\tilde{Z}_{q}=\sum_{k=1}^{q}\sum_{v\in \mathbb{D}_k}1_{\{S_v\geq 2Mq\}}.$$
Recall the definition for $M$ in (\ref{MaxMin}).
Let $Q'$ be defined by $\frac{dQ'}{dP}=e^{\lambda_+X_e-\Lambda(\lambda_+)}$. We have
\begin{eqnarray*}
  E\tilde{Z}_{q} & = & \sum_{k=1}^{q} b^k E1_{\{S_k\geq 2Mq\}}
  = \sum_{k=1}^{q} b^k E_{Q'} e^{-\lambda_+S_k+k\Lambda(\lambda_+)}1_{\{S_k\geq 2Mq\}}\\
  &\leq & \sum_{k=1}^{q} b^k  e^{-2\lambda_+Mq+k\Lambda(\lambda_+)} E_{Q'}1_{\{S_k\geq 2Mq\}}\\
  & \leq & \sum_{k=1}^{q} b^k e^{-\lambda_+Mk+k\Lambda(\lambda_+)}e^{-\lambda_+Mq}
  =e^{-\lambda_+Mq+o(q)},
\end{eqnarray*}
where, in the last equality, we use the definitions of $M$ and $\lambda_+$ (see (\ref{MaxMinCondition}) and (\ref{MaxMin})). It follows that
$$P(\tilde{Z}_{q}\geq 1)\leq E \tilde{Z}_{q}\leq e^{-\lambda_+Mq+o(q)}.$$
Again by the Borel-Cantelli lemma, $\tilde{Z}_{q}=0$ for all large $q$
almost surely. Taking
$q=\lfloor \epsilon_3 n^{1/3}\rfloor$ and combining with (\ref{LVUB}), we obtain that
$$L_n\leq L_{n+\lfloor
\epsilon_3n^{1/3}\rfloor}\leq (l_2+2M\epsilon_3) n^{1/3}$$
is true for all large $n$ almost surely. That is,
$$\limsup_{n\to \infty}\frac{L_n}{n^{1/3}}\leq l_2+2M\epsilon_3~~~~~~a.s..$$
Since $\epsilon_3>0$ and $l_2>l_0$ are arbitrary, we conclude that
\begin{equation}\label{ub}
  \limsup_{n\to \infty}\frac{L_n}{n^{1/3}}\leq l_0~~~~~~a.s..
\end{equation}
Together with (\ref{lb}), this completes the proof of Theorem
\ref{theo-main}.
\qed
\section{Concluding Remarks}

\subsection{The Curve s(t) of (\ref{DE})}

We comment in
this subsection on the appearance of the curve $s(t)$ of (\ref{DE}) as a solution to an appropriate variational principle. By the computation in Section 2, $s(t)n^{1/3}$ denotes the minimal possible position for vertices at level $tn$. However, in Section 3, it is not apriori clear that $s(t)$ will be our best choice. To see why this must indeed be
the best choice for the upper bound argument, let us consider a general curve $\phi(t)\leq l_2$ as the lower bound for the region. Examining the second
moment computation, we need
$$\max_{t}\{-\phi(t)+\int_0^t\frac{c}{(l_2-\phi(u))^2}du\}\leq 0$$
to make the argument work, where $c$ is some constant. Define $w(t)=l_2-\phi(t)\geq 0$. The above condition is equivalent to
$$l_2\geq \max_{t}\{w(t)+\int_0^t\frac{c}{w(u)^2}du\}.$$
Therefore, the best (smallest) upper bound that we can hope
is the result of the following optimization problem
\begin{equation}\label{OP}
\min_{w:(0,1)\to\mathds{R}_+}\max_{t}\{w(t)+\int_0^t\frac{c}{w(u)^2}du\}.
\end{equation}
The solution to this variational problem, denoted by $w^*(\cdot)$,
satisfies $s(t)=l_2-w^*(t)$.
\subsection{Generalizations}
The approach in this note seems to apply, under natural assumptions, to the situation where the $b$-ary tree is replaced by a Galton-Watson tree whose offspring distribution possesses high enough exponential moments. We do not pursue such an extension here.


\end{document}